\documentclass[10pt]{article}
\usepackage{amsmath,amssymb,enumerate,bbm}

\catcode`\à=\active \defà{\`a} \catcode`\À=\active \defÀ{\`A}
\catcode`\á=\active \defá{\'a} \catcode`\Á=\active \defÁ{\'A}
\catcode`\ä=\active \defä{\"a} \catcode`\Ä=\active \defÄ{\"A}
\catcode`\â=\active \defâ{\^a} \catcode`\Â=\active \defÂ{\^A}
\catcode`\å=\active \defå{{\aa}} \catcode`\Å=\active \defÅ{{\AA}}
\catcode`\ç=\active \defç{\c{c}} \catcode`\Ç=\active \defÇ{\c{C}}
\catcode`\è=\active \defè{\`e} \catcode`\È=\active \defÈ{\`E}
\catcode`\é=\active \defé{\'e} \catcode`\É=\active \defÉ{\'E}
\catcode`\ë=\active \defë{\"e} \catcode`\Ë=\active \defË{\"E}
\catcode`\ê=\active \defê{\^e} \catcode`\Ê=\active \defÊ{\^E}
\catcode`\ì=\active \defì{\`{\i}} \catcode`\Ì=\active \defÌ{\`{\I}}
\catcode`\í=\active \defí{\'{\i}} \catcode`\Í=\active \defÍ{\'{\I}}
\catcode`\ï=\active \defï{\"{\i}} \catcode`\Ï=\active \defÏ{\"{\I}}
\catcode`\î=\active \defî{\^{\i}} \catcode`\Î=\active \defÎ{\^{\I}}
\catcode`\ò=\active \defò{\`o} \catcode`\Ò=\active \defÒ{\`O}
\catcode`\ó=\active \defó{\'o} \catcode`\Ó=\active \defÓ{\'O}
\catcode`\ö=\active \defö{\"o} \catcode`\Ö=\active \defÖ{\"O}
\catcode`\ô=\active \defô{\^o} \catcode`\Ô=\active \defÔ{\^O}
\catcode`\ù=\active \defù{\`u} \catcode`\Ù=\active \defÙ{\`U}
\catcode`\ú=\active \defú{\'u} \catcode`\Ú=\active \defÚ{\'U}
\catcode`\ü=\active \defü{\"u} \catcode`\Ü=\active \defÜ{\"U}
\catcode`\û=\active \defû{\^u} \catcode`\Û=\active \defÛ{\^U}
\catcode`\ý=\active \defý{\'y} \catcode`\Ý=\active \defÝ{\'Y}
\catcode`\ÿ=\active \defÿ{\"y} \catcode`\˜=\active \def˜{\"Y}
\catcode`\½=\active \def½{!`}
\catcode`\¾=\active \def¾{?`}
\catcode`\ß=\active \defß{{\ss}}
\newcommand{\tmname}[1]{\textsc{#1}}
\newcommand{\tmem}[1]{{\em #1\/}}
\newcommand{\tmop}[1]{\ensuremath{\operatorname{#1}}}
\newtheorem{theorem}{Theorem}
\newtheorem{definition}{Definition}
\newcommand{\tmdummy}{$\mbox{}$}
\newenvironment{enumeratealpha}{\begin{enumerate}[a{\textup{)}}]}{\end{enumerate}}
\newcommand{\tmstrong}[1]{\textbf{#1}}
\newenvironment{enumerateroman}{\begin{enumerate}[i.]}{\end{enumerate}}
\newcommand{\nin}{\not\in}
\newtheorem{lemma}{Lemma}
\newenvironment{proof}{
  \noindent\textbf{Proof}\ }{\hspace*{\fill}
  \begin{math}\Box\end{math}\medskip}
\newenvironment{itemizeminus}
  {\begin{itemize}}{\end{itemize}}

\begin{document}

\title{{\tmname{Turing}} Computations on Ordinals} \author{Peter Koepke\\
University of Bonn} \maketitle

\begin{abstract}
  We define the notion of {\tmem{ordinal computability}} by generalizing
  standard {\tmname{Turing}} computability on tapes of length $\omega$ to
  computations on tapes of arbitrary ordinal length. We show that a set of
  ordinals is ordinal computable from a finite set of ordinal parameters if
  and only if it is an element of {\tmname{Gödel}}'s constructible universe
  $L$. This characterization can be used to prove the generalized continuum
  hypothesis in $L$.
\end{abstract}

\section{Introduction.}

A standard {\tmname{Turing}} computation may be visualized as a time-like
sequence of elementary {\tmem{read-write-move}} operations carried out by one
or more ``heads'' on ``tapes''. The sequence of actions is determined by the
initial tape contents and by a finite {\tmname{Turing}} {\tmem{program}}. The
specific choice of alphabet, operations and tapes may influence the time or
space complexity of calculations; by the {\tmname{Church-Turing}} thesis,
however, the associated notion of {\tmname{Turing}} {\tmem{computability}} is
not affected. So we may assume that {\tmname{Turing}} machines act on tapes
whose cells are indexed by the set $\omega$ ($=\mathbbm{N}$) of {\tmem{natural
numbers}} $0, 1, \ldots$ and contain $0$'s or $1$'s.

\begin{center}
  \begin{tabular}{|c|c|c|c|c|c|c|c|c|c|c|c|}
    \hline
    &  &  & S & P & A & C & E &  &  &  & \\
    \hline
    &  & 0 & 1 & 2 & 3 & 4 & 5 & 6 & 7 & $\ldots$ & $\ldots$ \\
    \hline
    & 0 & \underline{1} & 0 & 0 & 1 & 1 & 1 & 0 & 0 & 0 & 0\\
    \hline
    & 1 & 0 & \underline{0} & 0 & 1 & 1 & 1 & 0 & 0 &  & \\
    \hline
    T & 2 & 0 & 0 & \underline{0} & 1 & 1 & 1 & 0 & 0 &  & \\
    \hline
    I & 3 & 0 & \underline{0} & 1 & 1 & 1 & 1 & 0 & 0 &  & \\
    \hline
    M & 4 & 0 & 1 & \underline{1} & 1 & 1 & 1 & 0 & 0 &  & \\
    \hline
    E & : &  &  &  &  &  &  &  &  &  & \\
    \hline
    & $n$ & 1 & 1 & 1 & 1 & \underline{0} & 1 & 1 & 1 &  & \\
    \hline
    & $n + 1$ & 1 & 1 & 1 & 1 & 1 & \underline{1} & 1 & 1 &  & \\
    \hline
    & $\vdots$ &  &  &  &  &  &  &  &  &  & \\
    \hline
  \end{tabular}
\end{center}

\begin{center}
  {\small {\tmem{A standard {\tmname{Turing}} computation. Head positions are
  indicated by \underline{underlining}.}}}
\end{center}

An obvious generalization from the perspective of transfinite ordinal theory
is to extend {\tmname{Turing}} calculations to tapes whose cells are indexed
by the class $\tmop{Ord}$ of all {\tmem{ordinal numbers}}. Calculations will
become (infinite) sequences of elementary tape operations indexed by ordinals
which may be viewed as instances of time. For {\tmem{successor}} ordinals (or
times) calculations will basically be defined as for standard
{\tmname{Turing}} machines. At {\tmem{limit}} ordinals we define the tape
contents, program states and head positions by appropriate limit operations
which may be viewed as {\tmem{inferior limits}}.

\begin{center}
  \begin{tabular}{|c|c|c|c|c|c|c|c|c|c|c|c|c|c|c|c|}
    \hline
    &  & O & r & d & i & n & a & l &  & S & p & a & c & e & $\ldots$ \\
    \hline
    &  & 0 & 1 & 2 & 3 & 4 & 5 & 6 & 7 & $\ldots$ & $\ldots$  & $\omega$ &
    $\ldots$ & $\alpha$ & $\ldots$ \\
    \hline
    O & 0 & \underline{1} & 1 & 0 & 1 & 0 & 0 & 1 & 1 & $\ldots$ & $\ldots$ & 1 & $\ldots$
    & 1 & 0\\
    \hline
    r & 1 & 0 & \underline{1} & 0 & 1 & 0 & 0 & 1 & 1 &  &  & 1 &  &  & \\
    \hline
    d & 2 & 0 & 0 & \underline{0} & 1 & 0 & 0 & 1 & 1 &  &  & 1 &  &  & \\
    \hline
    i & 3 & 0 & \underline{0} & 0 & 1 & 0 & 0 & 1 & 1 &  &  & 1 &  &  & \\
    \hline
    n & 4 & 0 & 0 & \underline{0} & 0 & 0 & 0 & 1 & 1 &  &  & 1 &  &  & \\
    \hline
    a & : &  &  &  &  &  &  &  &  &  &  &  &  &  & \\
    \hline
    l & n & 1 & 1 & 1 & 1 & 0 & 1 & 0 & \underline{1} &  &  & 1 &  &  & \\
    \hline
    & n+1 & 1 & 1 & 1 & 1 & 1 & 1 & \underline{0} & 1 &  &  & 1 &  &  & \\
    \hline
    T & $\vdots$ & $\vdots$ & $\vdots$ & $\vdots$ & $\vdots$ & $\vdots$ &  & 
    &  &  &  &  &  &  & \\
    \hline
    i & $\omega$ & 0 & 0 & 1 & 0 & 0 & 0 & 1 & 1 & $\ldots$ & $\ldots$ & \underline{1} & 
    &  & \\
    \hline
    m & $\omega + 1$ & \underline{0} & 0 & 1 & 0 & 0 & 0 & 1 & 1 &  &  & 0 &  &  & \\
    \hline
    e & : &  &  &  &  &  &  &  &  &  &  &  &  &  & \\
    \hline
    $\vdots$ & $\theta$ & 1 & 0 & 0 & 1 & 1 & 1 & 1 & 0 & $\ldots$ & $\ldots$
    & $\ldots$ & $\ldots$ & \underline{0} & $\ldots$ \\
    \hline
    $\vdots$ & $\vdots$ &  &  & $\vdots$ &  &  & $\vdots$ &  &  & $\vdots$ &
    $\vdots$ &  &  &  & \\
    \hline
    $\vdots$ & $\vdots$ &  &  &  &  &  &  &  &  &  &  &  &  &  & \\
    \hline
  \end{tabular}
\end{center}

\begin{center}
  {\small {\tmem{An ordinal computation.}}}
\end{center}

The corresponding notion of {\tmem{ordinal computability}} obviously extends
{\tmname{Turing}} computability. By the {\tmname{Church-Turing}} thesis many
operations on natural numbers are ordinal computable. The {\tmem{ordinal
arithmetical}} operations (addition, multiplication, exponentiation) and other
basic operations on ordinals are also ordinal computable.

Indeed, the recursive properties of the family of ordinal computable functions
are so strong that the {\tmem{bounded truth predicate}}
\[ \{ ( \alpha, \varphi, \vec{x} ) | \alpha \in \tmop{Ord}, \varphi \text{ an
   } \in \text{-formula}, \vec{x} \in L_{\alpha}, L_{\alpha} \vDash \varphi (
   \vec{x} ) \} \]
for {\tmname{Gödel}}'s {\tmem{constructible hierarchy}} $L = \bigcup_{\alpha
\in \tmop{Ord}} L_{\alpha}$ is ordinal computable given some appropriate
coding. As a corollary we obtain the main result characterizing ordinal
computability:

\begin{theorem}
  A set $x \subseteq \tmop{Ord}$ is ordinal computable from finitely many
  ordinal parameters if and only if $x \in L$.
\end{theorem}

The implication from left to right will be immediate from the set-theoretical
absoluteness of ordinal computations. The converse requires a careful analysis
of the iterative definition of the constructible hierarchy to show that the
iteration can be carried out by an ordinal {\tmname{Turing}} machine.

This theorem may be viewed as an analogue of the {\tmname{Church-Turing}}
thesis: ordinal computability defines a natural and absolute class of sets,
and it is stable with respect to technical variations in its definition.

Theories of transfinite computations which extend {\tmname{Turing}}
computability have been proposed and studied for some time. {\tmem{Higher
recursion theory}} as described in the monograph {\cite{Sacks}} of
{\tmname{Gerald Sacks}} culminates in $E$-{\tmem{recursion}} which defines a
computational result $\{ e \} ( x )$ for programs $e$ built from basic set
functions applied to arbitrary sets $x$. The relation of $E$-computability to
constructibility is analogous to the statement of Theorem 1. In computer
science various infinitary machines like {\tmname{Büchi}} automata
{\cite{Buechi}} have been defined yielding important applications. The novelty
here is in keeping with the original {\tmname{Turing}} idea of reading and
writing on an unrestricted tape while interpreting unrestrictedness in the
widest sense as set-theoretical unboundedness.

Our work was inspired by the {\tmem{infinite time}} {\tmname{Turing}} machines
introduced by {\tmname{Joel D. Hamkins}}, {\tmname{Jeff Kidder}} and
{\tmname{Andy Lewis}} {\cite{HamkinsLewis}}. Infinite time {\tmname{Turing}}
machines use standard tapes indexed by natural numbers but allow infinite
computation sequences. At limit times, tape contents are defined as inferior
limits of previous contents. Inputs and outputs are characteristic functions
on the set $\omega$ of all natural numbers and may thus be viewed as
{\tmem{real numbers}}. The theory of infinite time {\tmname{Turing}} machines
is naturally related to definability theory over the structure $(\mathbbm{R},
\ldots )$, i.e., to {\tmem{descriptive set theory}}. In the case of tapes of
arbitrary ordinal length one is lead to consider a theory of arbitrarily
iterated definitions, i.e., {\tmem{constructibility theory}}.

\section{Ordinal {\tmname{Turing}} Machines}

We give an intuitive description of ordinal computations which will
subsequently be formalized. Consider a tape of ordertype $\tmop{Ord}$, i.e., a
sequence indexed by the class $\tmop{Ord}$ of all ordinals. The cells of the
tape can be identified with the ordinals, every cell can contain a $0$ or a
$1$ where $0$ is the default state. A {\tmem{read-write}} head moves on the
tape, starting at cell $0$. The computation is steered by a {\tmem{program}}
which consists of a finite sequence of {\tmem{commands}} indexed by natural
numbers. The indices of the commands can be seen as {\tmem{states}} of the
machine.

A computation of the machine is a sequence of machine configurations which are
indexed by ordinal ``times'' $0, 1, 2, \ldots, \omega, \omega + 1, \ldots$. At
time $t$ the read-write head reads the content of the cell at its position.
According to the content and the present machine state the head writes a
``$0$'' or a ``$1$'' and then moves to the right or to the left. Also the
machine changes into a new program state.

So far we have described the computation rules of finitary {\tmname{Turing}}
machines. {\tmem{Ordinal computations}} require the specification of the
behaviour at {\tmem{limit ordinals}}; we shall base the limit rules on simple
limit operations.

Assume that at time $t$ the head position is $H ( t )$. After a
{\tmem{move-right}} command we put $H ( t + 1 ) = H ( t ) + 1$. After a
{\tmem{move-left}} command we move one cell to the left if possible and
otherwise, if $H ( t )$ is a limit ordinal or $0$, jump to the default
position $0$:
\[ H ( t + 1 ) = \left\{\begin{array}{l}
     H ( t ) - 1, \text{ if } H ( t ) \text{ is a successor ordinal;}\\
     0, \text{ else.}
   \end{array}\right. \]
The definition of $H ( t )$ for $t$ a {\tmem{limit ordinal}} will be given
later.

At time $t$ the tape will be identified with a {\tmem{tape content}}
\[ T ( t ) = ( T ( t )_0, T ( t )_1, \ldots, T ( t )_{\omega}, T ( t )_{\omega
   + 1}, \ldots ) \]
which is a sequence of {\tmem{cell contents}} $T ( t )_{\alpha} \in \{ 0, 1
\}$. It is determined by previous write operations. For limit times $t$ the
content $T ( t )_{\alpha}$ of the $\alpha$-th cell is determined as follows:
if the cell content $T ( s )_{\alpha}$ stabilizes at a constant value as $s$
approaches $t$ we let $T ( t )_{\alpha}$ be that value; otherwise we take the
default value $T ( t )_{\alpha} = 0$. Formally this is an {\tmem{inferior
limit}}:
\[ T ( t )_{\alpha} = \liminf_{s \rightarrow t} T ( s )_{\alpha} . \]
A $\liminf$ rule will also be used for the program state and the head location
at limit times. Let $S ( t )$ be the program state at time $t$. For limit
times $t$ set
\[ S ( t ) = \liminf_{s \rightarrow t} S ( s ) . \]
Finally the head position $H ( t )$ for limit times $t$ is
\[ H ( t ) = \liminf_{s \rightarrow t, S ( s ) = S ( t )} H ( s ) . \]
The definitions of $S ( t )$ and $H ( t )$ can be motivated as follows. Since
a {\tmname{Turing}} program is finite its execution will lead to some
(complex) looping structure involving loops, subloops and so forth. This can
be presented by pseudo code like:

{\normalsize 
\begin{verbatim}
        ...
 17:begin loop
           ...
     21:   begin subloop
              ...
     29:   end subloop
           ...
     32:end loop
        ...
\end{verbatim}}

Assume that for times $s \rightarrow t$ the loop $( 17 - 32 )$ with its
subloop $( 21 - 29 )$ is traversed cofinally often. Then at limit time $t$ it
is natural to put the machine at the start of the ``main loop''. Assuming that
the lines of the program are enumerated in increasing order this corresponds
to the $\liminf$ rule
\[ S ( t ) = \liminf_{s \rightarrow t} S ( s ) . \]
The canonical head location $H ( t )$ is then determined as the inferior limit
of all head locations when the program is at the start of the ``main loop''.
If the head is for example moving linearly towards a limit location, say $H (
s_0 + i ) = h_0 + i$ for $i < \lambda$, we will have $H ( s_0 + \lambda ) =
h_0 + \lambda$. Note that the limit behaviour of the head position is defined
differently for infinite time {\tmname{Turing}} machines which do not possess
limit positions on the tape; there the head simply falls back to $0$ at limit
times.

The above intuitions are formalized as follows.

\begin{definition}
  \label{Def_Program} {\tmdummy}
  
  \begin{enumeratealpha}
    \item A {\tmstrong{command}} is a 5-tuple C=$( s, c, c', m, s' )$ where
    $s, s' \in \omega$ and $c, c', m \in \{ 0, 1 \}$; the natural number $s$
    is the {\tmstrong{state}} of the command $C$. The intention of the command
    $C$ is that if the machine is in state $s$ and reads the symbol $c$ under
    its read-write head, then it writes the symbol $c'$, moves the head left
    if $m = 0$ or right if $m = 1$, and goes into state $s'$. States
    correspond to the ``line numbers'' of some programming languages.
    
    \item A {\tmstrong{program}} is a finite set $P$ of commands satisfying
    the following structural conditions:
    \begin{enumerateroman}
      \item If $( s, c, c', m, s' ) \in P$ then there is $( s, d, d', n, t' )
      \in P$ with $c \neq d$; thus in state $s$ the machine can react to
      reading a ``$0$'' as well as to reading a ``$1$''.
      
      \item If $( s, c, c', m, s' ) \in P$ and $( s, c, c'', m', s'' ) \in P$
      then $c' = c'', m = m', s' = s''$; this means that the course of the
      computation is completely determined by the sequence of program states
      and the initial cell contents.
    \end{enumerateroman}
    \item For a program $P$ let
    \[ \tmop{states} ( P ) = \{ s | ( s, c, c', m, s' ) \in P \} \]
    be the set of program states.
  \end{enumeratealpha}
  
\end{definition}

\begin{definition}
  \label{ordinal_computation}Let $P$ be a program. A triple
  \[ S : \theta \rightarrow \omega, H : \theta \rightarrow \tmop{Ord}, T :
     \theta \rightarrow (^{\tmop{Ord}} 2 ) \]
  is an {\tmstrong{ordinal computation}} by $P$ if the following hold:
  \begin{enumeratealpha}
    \item $\theta$ is a successor ordinal or $\theta = \tmop{Ord}$; $\theta$
    is the {\tmstrong{length}} of the computation.
    
    \item $S ( 0 ) = H ( 0 ) = 0$; the machine starts in state $0$ with head
    position $0$.
    
    \item If $t < \theta$ and $S ( t ) \nin \tmop{state} ( P )$ then $\theta =
    t + 1$; the machine {\tmstrong{stops}} if the machine state is not a
    program state of $P$.
    
    \item If $t < \theta$ and $S ( t ) \in \tmop{state} ( P )$ then $t + 1 <
    \theta$; choose the unique command $( s, c, c', m, s' ) \in P$ with $S ( t
    ) = s$ and $T ( t )_{H ( t )} = c$; this command is executed as follows:
    \begin{eqnarray*}
      T ( t + 1 )_{\xi} & = & \left\{\begin{array}{l}
        c'  \text{, if } \xi = H ( t ) ;\\
        T ( t )_{\xi} \text{ , else;}
      \end{array}\right.\\
      S ( t + 1 ) & = & s' ;\\
      H ( t + 1 ) & = & \left\{\begin{array}{l}
        H ( t ) + 1 \text{, if } m = 1 ;\\
        H ( t ) - 1 \text{, if } m = 0 \text{ and } H ( t ) \text{ is a
        successor ordinal;}\\
        0 \text{, else.}
      \end{array}\right.
    \end{eqnarray*}
    \item If $t < \theta$ is a limit ordinal, the machine constellation at $t$
    is determined by taking inferior limits:
    \begin{eqnarray*}
      \forall \xi \in \tmop{Ord} T ( t )_{\xi} & = & \liminf_{r \rightarrow t}
      T ( r )_{\xi} ;\\
      S ( t ) & = & \liminf_{r \rightarrow t} S ( r ) ;\\
      H ( t ) & = & \liminf_{s \rightarrow t, S ( s ) = S ( t )} H ( s ) .
    \end{eqnarray*}
  \end{enumeratealpha}
  The computation is obviously recursively determined by the initial tape
  contents $T ( 0 )$ and the program $P$. We call it the {\tmstrong{ordinal
  computation}} {\tmstrong{by}} $P$ {\tmstrong{with input}} $T ( 0 )$. If the
  computation stops, $\theta = \beta + 1$ is a successor ordinal and $T (
  \beta )$ is the final tape content. In this case we say that $P$
  {\tmstrong{computes}} $T ( \beta )$ from $T ( 0 )$ and write $P : T ( 0 )
  \mapsto T ( \beta )$.
\end{definition}

This interpretation of programs yields associated notions of computability.

\begin{definition}
  A partial function $F :^{\tmop{Ord}} 2 \rightharpoonup^{\tmop{Ord}} 2$ is
  {\tmstrong{ordinal computable}} if there is a program $P$ such that $P : T
  \mapsto F ( T )$ for every $T \in \tmop{dom} ( F )$. 
\end{definition}

By coding, the notion of ordinal computability can be extended to other
domains. We can e.g. {\tmem{code}} an ordinal $\delta \in \tmop{Ord}$ by the
characteristic function $\chi_{\{ \delta \}} : \tmop{Ord} \rightarrow 2$,
$\chi_{\{ \delta \}} ( \xi ) = 1$ iff $\xi = \delta$, and define:

\begin{definition}
  A partial function $F : \tmop{Ord} \rightharpoonup \tmop{Ord}$ is
  {\tmstrong{ordinal computable}} if the function $\chi_{\{ \delta \}} \mapsto
  \chi_{\{ F ( \delta ) \}}$ is ordinal computable.
\end{definition}

We also consider computations involving finitely many ordinal
{\tmem{parameters}}.

\begin{definition}
  A subset $x \subseteq \tmop{Ord}$ is {\tmstrong{ordinal computable from
  finitely many ordinal parameters}} if there a finite subset $z \subseteq
  \tmop{Ord}$ and a program $P$ such that $P : \chi_z \mapsto \chi_x$.
\end{definition}

In view of our intended applications of ordinal computations to models of set
theory we note some absoluteness properties:

\begin{lemma}
  \label{absolutness}Let $( M, \in )$ be a transitive model of $\tmop{ZF}^-$,
  i.e., of {\tmname{Zermelo-Fraen{\-}kel}} set theory without the powerset axiom.
  Let $P$ be a program and let $T ( 0 ) : \tmop{Ord} \rightarrow 2$ be an
  initial tape content so that $T ( 0 ) \upharpoonright ( \tmop{Ord} \cap M )$
  is definable in $M$. Let $S : \theta \rightarrow \omega, H : \theta
  \rightarrow \tmop{Ord}, T : \theta \rightarrow (^{\tmop{Ord}} 2 )$ be the
  ordinal computation by $P$ with input $T ( 0 )$. Then:
  \begin{enumeratealpha}
    \item The ordinal computation by $P$ with input $T ( 0 )$ is absolute for
    $( M, \in )$ below $( \tmop{Ord} \cap M )$, i.e.,
    \[ S : \theta \cap M \rightarrow \omega, H : \theta \cap M \rightarrow
       \tmop{Ord}, \bar{T} : \theta \cap M \rightarrow (^{\tmop{Ord} \cap M} 2
       ) \]
    with $\bar{T} ( t ) = T ( t ) \upharpoonright ( \tmop{Ord} \cap M )$ is
    the ordinal computation by $P$ with input $T ( 0 ) \upharpoonright (
    \tmop{Ord} \cap M )$ as computed in the model $( M, \in )$.
    
    \item If $\tmop{Ord} \subseteq M$ then the ordinal computations by $P$ in
    $M$ and in the universe $V$ are equal.
    
    \item Let $\tmop{Ord} \subseteq M$ and $x, y \subseteq \tmop{Ord}$, $x, y
    \in M$. Then $P : \chi_x \mapsto \chi_y$ if and only if $( M, \in )
    \vDash$``$P : \chi_x \mapsto \chi_y$''.
    
    \item Let $x, y \subseteq \tmop{Ord}$, $x, y \in M$. Assume that $( M,
    \in ) \vDash$``$P : \chi_x \mapsto \chi_y$''. Then $P : \chi_x \mapsto
    \chi_y $.
  \end{enumeratealpha}
\end{lemma}

The properties follow from the observation that the recursion in Definition
\ref{ordinal_computation} is clearly absolute between $M$ and $V$. Note that
the converse of d) is in general false. With the subsequent results on
constructibility we could let $M = L_{\delta}$ be the minimal level of the
constructible hierarchy which is a model of $\tmop{ZF}^-$. If $P$ is a program
which searches for the minimal ordinal $\delta$ such that $L_{\delta}$ is a
$\tmop{ZF}^-$-model then $P$ will stop in $V$ but not in $M$.

\section{Ordinal Algorithms}

We present a number of fundamental algorithms which can be implemented as
ordinal computations. Our emphasis is not on writing concrete
{\tmem{programs}} as in Definition \ref{Def_Program} but on showing that
programs {\tmem{exist}}. It thus suffices to present basic ideas and
algorithms together with methods to combine these into complex algorithms. We
shall freely use informal ``higher programming languages'' to describe
algorithms. Algorithms are based on data formats for the representation of
input and output values. Again we shall not give detailed definitions but only
indicate crucial features of the formats.

The intended computations will deal with ordinals and sequences of ordinals.
The simplest way of representing the ordinal $\alpha \in \tmop{Ord}$ in an
ordinal machine is by a tape whose content is the characteristic function of
$\{ \alpha \}$:
\[ \text{$\chi_{\{ \alpha \}} : \tmop{Ord} \rightarrow 2$,  $\chi_{\{ \alpha
   \}} ( \xi ) = 1$ iff $\xi = \alpha$.} \]
A basic task is to {\tmem{find}} or {\tmem{identify}} this ordinal $\alpha$:
initially the head is in position $0$, it then moves to the right until it
stops exactly at position $\alpha$. This is achieved by the following program:
\[ P = \{ ( 0, 0, 0, 1, 0 ), ( 0, 1, 1, 1, 1 ), ( 1, 0, 0, 0, 2 ), ( 1, 1, 1,
   0, 2 ) \} . \]
The program is in state $0$ until it reads a $1$, then it goes one cell to the
right, one cell to the left, and stops because $2$ is not a program state.
Informally the algorithm may be written as
\begin{verbatim}
Find_Ordinal:
   if head = 1 then STOP otherwise moveright
\end{verbatim}
Similarly one can let the head find (the beginning) of any finite
0-1-{\tmem{bitstring}} $b_0 \ldots b_{k - 1} $:
\begin{verbatim}
Find_Bitstring :
A:     if head = `' then moveright otherwise goto C0
       if head = `' then moveright otherwise goto C1
       ...
       if head = `' then goto B otherwise goto C(k-1)
B:     moveleft
       ...
       moveleft
       moveleft
       stop
C(k-1):moveleft
       ...
C1:    moveleft
C0:    moveright
       goto A
\end{verbatim}
In view of this algorithm we may assume that the tape contains arbitrary
{\tmem{symbols}} coded by finite bitstrings instead of single bits. Note that
the above programs obviously perform the intended tasks on standard
{\tmname{Turing}} machines. The limit rules are designed to lift this
behaviour continuously to transfinite ordinals.

Often one has to reset the head to its initial position $0$. There are several
methods to achieve this. A universal one assumes that there is a unique
initial inscription \verb¤start¤ on the tape which indicates the $0$-position:
\begin{verbatim}
Reset_head:
A:     moveleft
       if head reads `start' then STOP otherwise goto A 
\end{verbatim}
It will be convenient to work with several tapes side-by-side instead of just
one. This corresponds to the idea of program {\tmem{variables}} whose values
are checked and manipulated. One can simulate an $n$-tape machine on a
$1$-tape machine. The contents $( T^i_{\xi} | \xi \in \tmop{Ord} )$ of the
$i$-th tape are successively written into the cells of tape $T$ indexed by
ordinals $2 n \xi + 2 i$:
\[ T_{2 n \xi + 2 i} = T^i_{\xi} . \]
The head position $H^i$ on the $i$-th tape is simulated by writing 1's into an
initial segment of length $H^i$ of cells with indices of the form $2 n \xi + 2
i + 1$:
\[ T_{2 n \xi + 2 i + 1} = \left\{\begin{array}{l}
     1 \text{, if } \xi < H^i ;\\
     0 \text{, else} .
   \end{array}\right. \]
So two tapes with contents $a_0 a_1 a_2 a_3 a_4 \ldots$ and $b_0 b_1 b_2 b_3
b_4 \ldots$ and head positions $3$ and 1 respectively are coded as
\[ T = a_0 1 b_0 1 a_1 1 b_1 0 a_2 1 b_2 0 a_3 0 b_3 0 a_4 0 b_4 0 \ldots
   \ldots . \]
We describe operations of machines with several tapes by commands like
\verb¤move-¤ \verb¤right2¤ or \verb¤print3 = ` '¤, where the number of the active tape
is adjoined to the right. There are canonical but tedious translations from
programs for $n$-tape machines into corresponding programs for $1$-tape
machines. A manipulation of the $i$-th tape amounts to first finding the head
marker at ordinals of form $2 n \xi + 2 i + 1$; moving left by one cell one
obtains the corresponding cell content for possible modification; the
subsequent head movement is simulated by moving right again, writing a $0$,
moving $2 n$ cells to the right or left, and printing a $1$; if a
left-movement goes across a limit ordinal, then a ``$1$'' has to be printed
into cell $2 i + 1$.

The subsequent algorithms will be presented as multiple tape algorithms. One
can assume that one or more of the tapes serve as standard {\tmname{Turing}}
tapes on which ordinary {\tmname{Turing}} recursive functions are computed.
Since the usual syntactical operations for a language of set theory are
intuitively computable we can assume by the {\tmname{Church-Turing}} thesis
that these operations are ordinal computed on some of the ordinal tapes. This
will be used in the ordinal computation of the constructible model $L$.

Basic operations on ordinals are ordinal computable. Let the ordinals $\alpha$
and $\beta$ be given on tapes $0$ and $1$ as their characteristic functions
$\chi_{\{ \alpha \}}$ and $\chi_{\{ \beta \}} $. The following algorithm
{\tmem{compares}} the ordinals and indicates the result of the comparison by
its ``stopping state'':
\begin{verbatim}
Ordinal_Comparison:
       Reset_Head0
       Reset_Head1
A:     if head0 = `1' and head1 = `0' then STOP (`alpha > beta')
       if head0 = `1' and head1 = `1' then STOP (`alpha = beta')
       if head0 = `0' and head1 = `1' then STOP (`alpha < beta')
       moveright0
       moveright1
       goto A
\end{verbatim}
Obviously there are ordinal algorithms to {\tmem{reset}} a register containing
an ordinal to $0$, or to {\tmem{copy}} one ordinal register to another one.
The ordinal {\tmem{sum}} $\alpha + \beta$ and {\tmem{product}} $\alpha \cdot
\beta$ are computable as follows:
\begin{verbatim}
Ordinal_Addition:
       Reset_Head0
       Reset_Head1
       Reset_Head2
A:     if head0 = `1' then goto B
       moveright0
       moveright2
       goto A
B:     if head1 = `1' then goto C
       moveright1
       moveright2
       goto B
C:     print2 = `1'
       STOP
\end{verbatim}
\begin{verbatim}
Ordinal_Multiplication:
       Reset_Head0
       Reset_Head1
       Reset_Head2
A:     if head1 = `1' then goto C
       if head0 = `1' then goto B
       moveright0
       moveright2
       goto A
B:     Reset_Head0
       moveright1
       goto A
C:     print2 = `1'
       STOP
\end{verbatim}

The class $\tmop{Ord}^{< \omega} = \{ s | \exists k < \omega s : k
\rightarrow \tmop{Ord} \}$ of finite sequences of ordinals will be of
particular interest for relating ordinal computability to the iterated
definability of {\tmname{Gödel}}'s constructible universe. We code a sequence
$( \alpha_0, \ldots, \alpha_{k - 1} ) : k \rightarrow \tmop{Ord}$ by a tape
which starts with an initial symbol ``('', followed by $k$ intervals of 0's of
lengths $\alpha_0, \ldots \alpha_{k - 1}$ respectively, which are separated by
a separation symbol ``,'' and then a closing ``)''. So $( 1, \omega, \omega +
2 )$ is coded as
\[ ( 0, 00 \ldots, 00 \ldots 00 ) \]
If the sequence is given on tape 0 and a natural number $n$ on tape 1 then the
$n$-th element of the sequence can be output on tape 2 by the following
algorithm:
\begin{verbatim}
Extract:
       Reset_Head0
       Reset_Head1
       Reset_Head2
A:     if head1 = `1' then goto C
       moveright1
B:     if head0 = `,' then goto A
       if head0 = `)' then STOP (no output)
       moveright0
       goto B
C:     moveright0
D:     if head0 = `,' then goto E
       if head0 = `)' then goto E
       moveright0
       moveright2
       goto D
E:     print2 = `1'
       STOP
\end{verbatim}

Another important operation on sequences is the {\tmem{replacement}} of the
$n$-th element of a sequence $s$ of ordinals by a given ordinal $\alpha$; if
the given sequence is shorter than $n + 1$, it is padded by 0's up to length
$n + 1$. Formally this operation on sequences is defined as $s \mapsto s
\frac{\alpha}{n}$ where $\tmop{dom} ( s \frac{\alpha}{n} ) = \tmop{dom} ( s )
\cup ( n + 1 )$ and
\[ s \frac{\alpha}{n} ( i ) = \left\{\begin{array}{l}
     s ( i ) \text{, if } i \in \tmop{dom} ( s ) \setminus \{ n \} ;\\
     \alpha \text{, if } i = n ;\\
     0 \text{, else} .
   \end{array}\right. \]
Let the original sequence be given on tape 0, the natural number $n$ on tape
1, and the ordinal $\alpha$ on tape 2. The modified sequence $s
\frac{\alpha}{n}$ can be output on tape 3 by the following algorithm:
\begin{verbatim}
Replace:
       print3 = `('
A:     moveright0
       if read1 = `1' then goto C
       if read0 = `,' then goto B
       if read0 = `)' then goto H
       moveright3
       goto A
B:     print3 = `,'
       moveright1
       goto A
C:     if read2 = `1' then goto D
       moveright2
       moveright3
       goto C
D:     moveright0
       if read0 = `0' then goto D
E:     if read0 = `,' then goto F
       if read0 = `)' then goto G
       moveright3
       moveright0
       goto E
F:     print3 = `,'
       moveright3
       moveright0
       goto E
G:     print3 = `)'
       STOP
H:     print3 = `,'
       moveright1
       if read1 = `0' then goto H
I:     if read2 =`1' then goto J
       moveright2
       moveright3
       goto I
J:     print3 = `)'
       STOP
\end{verbatim}
With the subroutine mechanism known from ordinary programming the basic
algorithms can be combined into complex algorithms for comparing and
manipulating ordinal sequences. We can, e.g., carry out a syntactic
manipulation on a standard {\tmname{Turing}} tape which outputs requests for
checking or manipulating elements of ordinal sequences. According to the
requests the appropriate elements can be extracted and subjected to some
algorithms whose results can be substituted into the original sequences.

\section{Enumerating Finite Sequences of Ordinals}

For $X$ a class let $[ X ]^{< \omega} = \{ z \subseteq X | z \text{ is finite}
\}$ and  $X^{< \omega} = \{ s | \exists k < \omega s : k \rightarrow X \}$ be
the class of all finite {\tmem{subsets}} of $X$ and of all finite
{\tmem{sequences}} from $X$ respectively.

Finite sequences of ordinals are finite sets of ordered pairs:
\[ \text{$\tmop{Ord}^{< \omega} \subseteq [ \omega \times \tmop{Ord} ]^{<
   \omega}$} . \]
Well-order $\omega \times \tmop{Ord}$ by
\[ ( m, \alpha ) \prec ( n, \beta ) \text{ iff } \alpha < \beta \text{ or } (
   \alpha = \beta \wedge m < n ) . \]
Define a canonical well-order $( [ \omega \times \tmop{Ord} ]^{< \omega},
\prec^{\ast} )$ by largest difference:
\[ s \prec^{\ast} s' \text{  iff  } \exists x \in s' \setminus s \{ y \in s |
   y \succ x \} = \{ y \in s' | y \succ x \} . \]
One can show inductively that $( [ Y ]^{< \omega}, \prec^{\ast} )$ is a
well-order on initial segments $Y$ of $( \omega \times \tmop{Ord}, \prec )$.
So $\prec^{\ast}$ well-orders $[ \omega \times \tmop{Ord} ]^{< \omega}$ and
hence $\tmop{Ord}^{< \omega}$. We note an important {\tmem{substitution
property}} of the well-order:

\begin{lemma}
  \label{substition_property}If $s, t \in \tmop{Ord}^{< \omega}$, $m \in
  \tmop{dom} ( s )$, $\tmop{dom} ( s ) \subseteq \tmop{dom} ( t )$, $s
  \upharpoonright ( \tmop{dom} ( s ) \setminus \{ m \} ) = t \upharpoonright (
  \tmop{dom} ( s ) \setminus \{ m \} )$, $t ( m ) < s ( m )$, $\forall i \in
  \tmop{dom} ( t ) \setminus \tmop{dom} ( s ) t ( i ) < s ( m )$ then
  \[ t \prec^{\ast} s. \]
\end{lemma}

So replacing an arbitrary ordinal $s ( m )$ of $s$ by possibly many smaller
ordinals leads to a descent in $\prec^{\ast}$. The substitution property will
correspond to the substitution of a bounded variable below some bound by terms
with parameters smaller than that bound. This will lead to a recursive
definition of bounded truth in $L$ along the $\prec^{\ast}$-relation.

We define an enumeration $S : \tmop{Ord} \rightarrow \tmop{Ord}^{< \omega}$
(with repetitions) of $\tmop{Ord}^{< \omega}$ which is compatible with
$\prec^{\ast}$ and which can be computed by an ordinal machine. The idea of
the construction is to recursively apply the replacement operation $s
\frac{\alpha}{n}$ to sequences $s$ which have been enumerated before.

For $( m, \alpha ) \in \omega \times \tmop{Ord}$ define functions $S_{m
\alpha} : \theta_{m \alpha} \rightarrow \tmop{Ord}^{< \omega}$ such that for
$( m, \alpha ) \prec ( n, \beta )$, $S_{m \alpha}$ is an initial segment of
$S_{n \beta}$. Set $S_{0 0} : 1 \rightarrow \tmop{Ord}^{< \omega}$, $S_{0
0} ( 0 ) = \emptyset$. For $\beta > 0$ set
\[ S_{0 \beta} = \bigcup_{( m, \alpha ) \prec ( 0, \beta )} S_{m \alpha} ._{}
\]
Assume that $S_{m \alpha} : \theta_{m \alpha} \rightarrow \tmop{Ord}^{<
\omega}$ is defined. Then define $S_{m + 1, \alpha} : \theta_{m \alpha} \cdot
2 \rightarrow \tmop{Ord}^{< \omega}$ by: $S_{m + 1, \alpha} \upharpoonright
\theta_{m \alpha} = S_{m \alpha} $; for $\xi < \theta_{m \alpha}$ let
\[ S_{m + 1, \alpha} ( \theta_{m \alpha} + \xi ) = S_{m \alpha} ( \xi )
   \frac{\alpha}{m} . \]
Finally set
\[ S_{} = \bigcup_{( m, \alpha ) \in \omega \times \tmop{Ord}} S_{m \alpha}
   ._{} \]
\begin{lemma}
  \begin{enumeratealpha}
    \item $S : \tmop{Ord} \rightarrow \tmop{Ord}^{< \omega}$ is a surjection.
    
    \item If $\xi < \zeta$ then $S ( \xi ) = S ( \zeta )$ or  $S ( \xi )
    \prec^{\ast} S ( \zeta )$.
  \end{enumeratealpha}
\end{lemma}

\begin{proof}
  a) We show by induction on $\alpha$ that $S_{0 \alpha} : \theta_{0 \alpha}
  \rightarrow \alpha^{< \omega}$ is a surjection. The initial case $\alpha =
  0$ and the limit step are easy. Consider $\alpha = \beta + 1$ and some $s
  \in \alpha^{< \omega}$, $s : k \rightarrow \alpha$. Let $\bar{s} : k
  \rightarrow \alpha$ be the following restriction of $s$ to $\beta$:
  \[ \bar{s} ( i ) = \left\{\begin{array}{l}
       s ( i ), \text{ if } s ( i ) < \beta ;\\
       0, \text{ if } s ( i ) = \beta .
     \end{array}\right. \]
  By the inductive assumption there is $\xi < \theta_{0 \beta}$ such that
  $S_{0 \beta} ( \xi ) = \bar{s}$. Then
  \[ S_{0 \alpha} ( \theta_{0 \beta} \cdot ( \sum_{i < k, s ( i ) = \beta} 2^i
     ) + \xi ) = s . \]
  b) follows from the substitution property.
\end{proof}

The enumeration $S$ of $\tmop{Ord}^{< \omega}$ is {\tmem{ordinal computable}}
using coding methods from the previous paragraph. We indicate a program which
writes the values of $S$ consecutively on a tape:
\[ S ( 0 ) S ( 1 ) S ( 2 ) \ldots S ( \omega ) S ( \omega + 1 ) \ldots S (
   \omega + \omega ) \ldots S ( \alpha ) \ldots \ldots \]
where each $S ( \alpha )$ is of the form
\[ ( 0 \ldots 0, 0 \ldots 0, \ldots \ldots, 0 \ldots 0 ) \]
The algorithm is based on the \texttt{Replace}-algorithm from the previous
section:
\begin{verbatim}
counter0 = 0
counter1 = 0
position_of_writing_head = 0
write the empty sequence `()'
while true
 while counter1 < 
  mark = position_of_writing_head
  position_of_reading_head = 0
  while position_of_reading_head < mark
   read sequence
   Replace element at position counter1 by counter0 
   write modified sequence at mark
  endwhile
  counter1 = counter1 + 1
 endwhile
 counter1 = 0
 counter0 = counter0 + 1
endwhile
\end{verbatim}
The procedure will eventually be extended as to write a bounded truth function
for the constructible hierarchy.

\section{The Constructible Hierarchy}

{\tmname{Kurt Gödel}} {\cite{Goedel}} defined the inner model $L$ of
{\tmem{constructible sets}} as the union of a hierarchy of levels $L_{\alpha}
$:
\[ L = \bigcup_{\alpha \in \tmop{Ord}} L_{\alpha}  \]
where the hierarchy is defined by: $L_0 = \emptyset$, $L_{\delta} =
\bigcup_{\alpha < \delta} L_{\alpha}$ for limit ordinals $\delta$, and
$L_{\alpha + 1} =$the set of all sets which are first-order definable in the
structure $( L_{\alpha}, \in )$. The standard reference to the theory of the
model $L$ is the book {\cite{Devlin}} by {\tmname{Keith Devlin}}.

An element of $L$ is definable over some $L_{\alpha}$ from parameters which
are themselves definable over some $L_{\beta} $, $\beta < \alpha$ in some
other parameters and so forth. We therefore introduce a language with
definable terms, which in turn may involve definable terms etc.

Consider a language with symbols $(, ), \{, \}, |, \in, =, \wedge, \neg,
\forall, \exists$ and variables $v_0, v_1, \ldots$. We define
({\tmstrong{bounded}}) {\tmstrong{formulas}} and ({\tmstrong{bounded}})
{\tmstrong{terms}} by a common recursion on the lenghts of words formed from
these symbols:
\begin{itemizeminus}
  \item the variables $v_0, v_1, \ldots$ are terms;
  
  \item if $s$ and $t$ are terms then $s = t$ and $s \in t$ are formulas;
  
  \item if $\varphi$ and $\psi$ are formulas then $\neg \varphi$, $( \varphi
  \wedge \psi )$, $\forall v_i \in v_j \varphi$ and $\exists v_i \in v_j
  \varphi$ are formulas;
  
  \item if $\varphi$ is a formula then $\{ v_i \in v_j | \varphi \}$ is a
  term.
\end{itemizeminus}
For terms and formulas of this language define {\tmstrong{free}} and
{\tmstrong{bound}} {\tmstrong{variables}}:
\begin{itemizeminus}
  \item $\tmop{free} ( v_i ) = \{ v_i \}, \tmop{bound} ( v_i ) = \emptyset$;
  
  \item $\tmop{free} ( s = t ) = \tmop{free} ( s \in t ) = \tmop{free} ( s )
  \cup \tmop{free} ( t )$;
  
  \item $\tmop{bound} ( s = t ) = \tmop{bound} ( s \in t ) = \tmop{bound} ( s
  ) \cup \tmop{bound} ( t )$;
  
  \item $\tmop{free} ( \neg \varphi ) = \tmop{free} ( \varphi ), \tmop{bound}
  ( \neg \varphi ) = \tmop{bound} ( \varphi )$;
  
  \item $\tmop{free} ( ( \varphi \wedge \psi ) ) = \tmop{free} ( \varphi )
  \cup \tmop{free} ( \psi ), \tmop{bound} ( ( \varphi \wedge \psi ) ) =
  \tmop{bound} ( \varphi ) \cup \tmop{bound} ( \psi )$;
  
  \item $\tmop{free} ( \forall v_i \in v_j \varphi ) = \tmop{free} ( \exists
  v_i \in v_j \varphi ) = \tmop{free} ( \{ v_i \in v_j | \varphi \} ) = (
  \tmop{free} ( \varphi ) \cup \{ v_j \} ) \setminus \{ v_i \}$;
  
  \item $\tmop{bound} ( \forall v_i \in v_j \varphi ) = \tmop{bound} ( \exists
  v_i \in v_j \varphi ) = \tmop{bound} ( \{ v_i \in v_j | \varphi \} ) = \\
  = \tmop{bound} ( \varphi ) \cup \{ v_i \}$.
\end{itemizeminus}

For technical reasons we will be interested in terms and formulas in which
\begin{itemizeminus}
  \item no bound variable occurs free,
  
  \item every free variable occurs exactly once.
\end{itemizeminus}
Such terms and formulas are called {\tmstrong{tidy}}; with tidy formulas one
avoids having to deal with the interpretation of one free variable at
different positions within a formula.

In recursive truth definitions one reduces the truth of formulas to the truth
of {\tmem{simpler}} formulas. The {\tmstrong{term complexity}} $\tmop{tc} ( t
)$ and $\tmop{tc} ( \varphi )$ of terms and formulas is defined recursively:
\begin{itemizeminus}
  \item $\tmop{tc} ( v_i ) = 0$;
  
  \item $\tmop{tc} ( s = t ) = \tmop{tc} ( s \in t ) = \max ( \tmop{tc} ( s ),
  \tmop{tc} ( t ) )$;
  
  \item $\tmop{tc} ( \neg \varphi ) = \tmop{tc} ( \forall v_i \in v_j \varphi
  ) = \tmop{tc} ( \exists v_i \in v_j \varphi ) = \tmop{tc} ( \varphi )$;
  
  \item $\tmop{tc} ( \varphi \wedge \psi ) = \max ( \tmop{tc} ( \varphi ),
  \tmop{tc} ( \psi ) )$;
  
  \item $\tmop{tc} ( \{ v_i \in v_j | \varphi \} ) = \tmop{tc} ( \varphi ) +
  1$.
\end{itemizeminus}

We can define a pre-wellordering $<_{\tmop{Form}}$ of the set of all bounded
formulas by
\[ \varphi <_{\tmop{Form}} \psi \text{ iff } \tmop{tc} ( \varphi ) < \tmop{tc}
   ( \psi ) \text{ or } ( \tmop{tc} ( \varphi ) = \tmop{tc} ( \psi ) \wedge
   \tmop{length} ( \varphi ) < \tmop{length} ( \psi ) ) . \]
Obviously the syntactical notions and operations of this language are
{\tmname{Turing}} computable and therefore ordinal computable. Also there is
an ordinal computable enumeration of all formulas which is compatible with
$<_{\tmop{Form}} $.

An {\tmem{assignment}} for a term $t$ or formula $\varphi$ is a finite
sequence $a : k \rightarrow V$ so that for every free variable $v_i$ of $t$ or
$\varphi$ we have $i < k$; $a ( i )$ will be the {\tmem{interpretation}} of
$v_i $. The {\tmem{value}} of $t$ or the {\tmem{truth value}} of $\varphi$ is
determined by the assignment $a$. We write $t [ a ]$ and $\varphi [ a ]$ for
the values of $t$ und $\varphi$ under the assignment $a$.

Concerning the constructible hierarchy $L$, it is shown by an easy induction
on $\alpha$ that every element of $L_{\alpha}$ is the interpretation $t [ (
L_{\alpha_0}, L_{\alpha_1}, \ldots, L_{\alpha_{k - 1}} ) ]$ of some
{\tmem{tidy}} term $t$ with an assignment $( L_{\alpha_0}, L_{\alpha_1},
\ldots, L_{\alpha_{k - 1}} )$ whose values are constructible levels
$L_{\alpha_i}$ with $\alpha_0, \ldots, \alpha_{k - 1} < \alpha$. This will
allow to reduce bounded quantifications $\forall v \in L_{\alpha}$ or $\exists
v \in L_{\alpha}$ to the substitution of terms of lesser complexity. Moreover,
the truth of (bounded) formulas in $L$ is captured by {\tmem{tidy}} bounded
formulas of the form $\varphi [ ( L_{\alpha_0}, L_{\alpha_1}, \ldots,
L_{\alpha_{k - 1}} ) ]$.

We shall code an assignment of the form $( L_{\alpha_0}, L_{\alpha_1},
\ldots, L_{\alpha_{k - 1}} )$ by its sequence of ordinal indices, i.e., we
write 

\[t [ ( \alpha_0, \alpha_1, \ldots, \alpha_{k - 1} ) ] \text{ or } 
\varphi [ (
\alpha_0, \alpha_1, \ldots, \alpha_{k - 1} ) ]\] 
instead of 
\[t [ (
L_{\alpha_0}, L_{\alpha_1}, \ldots, L_{\alpha_{k - 1}} ) ]
\text{ or } \varphi [ (
L_{\alpha_0}, L_{\alpha_1}, \ldots, L_{\alpha_{k - 1}} ) ].\] 
The relevant
assignments are thus elements of $\tmop{Ord}^{< \omega}$ and can be handled by
the programs of the previous section. Since the bounded language is recursive
we can modify the enumeration program so that all assigned tidy formulas
$\varphi [ ( \alpha_0, \alpha_1, \ldots, \alpha_{k - 1} ) ]$ occur in the
enumeration: for a fixed assigment $a = ( \alpha_0, \alpha_1, \ldots,
\alpha_{k - 1} )$ list the pairs $( a, \varphi )$ where $\varphi$ is a tidy
formula with $\tmop{free} ( \varphi ) \subseteq k$ in an order compatible with
$<_{\tmop{Form}} $. The following is a straightforward extension of the
enumeration program of the previous section:
\begin{verbatim}
counter0 = 0
counter1 = 0
position_of_writing_head = 0
write the empty sequence `()'
while true
 while counter1 < 
  mark = position_of_writing_head
  position_of_reading_head = 0
  while position_of_reading_head < mark
   read sequence
   Replace element at position counter1 by counter0
   form an enumeration of the appropriate tidy formulas 
        which is compatible with the pre-wellorder   
   for all enumerated formulas
    write the modified sequence and the formula
   endfor 
  endwhile
  counter1 = counter1 + 1
 endwhile
 counter1 = 0
 counter0 = counter0 + 1
endwhile
\end{verbatim}

\section{A Bounded Truth Function for $L$}

We define a bounded truth function $W$ for the constructible hierarchy on the
class
\[ A = \text{$\{ ( a, \varphi ) | a \in \tmop{Ord}^{< \omega}, \varphi \text{
   is a tidy bounded formula}, \tmop{free} ( \varphi ) \subseteq \tmop{dom} (
   a )$} \} \]
of all ``tidy pairs'' of assignments and formulas. Define the
{\tmstrong{bounded constructible truth function}} $W : A \rightarrow 2$ by
\[ W ( a, \varphi ) = 1 \text{ iff } \varphi [ a ] . \]
The function $W$ has a {\tmem{recursive}} definition along the enumeration of
$A$ given by the enumeration algorithm from the preceding section. We explain
the principal idea of the recursion with a bounded quantification like
$\exists v_i \in v_j \varphi [ a ]$. If $a ( j ) = \alpha$ then the assigned
formula is satisfied if and only if there is a witness for $\varphi$ in
$L_{\alpha} $. By the recursive definition of $L_{\alpha}$ such a witness must
be the interpretation $t [ b ]$ of a term where $b ( l ) < \alpha$ for every
free variable $v_l$ of $t$. If one chooses $t$ such that it has no variable in
common with $\varphi$ the assignment $b$ can be taken to further satisfy $b (
j ) < a ( j ) = \alpha$. The substitution property of Lemma
\ref{substition_property} leads to the evaluation of $\varphi \frac{t}{v_i} [
b ]$ for some $b \prec^{\ast} a$ which is the basis for the subsequent
recursion.

As we want to work with {\tmem{tidy}} formulas a technical problem has to be
solved. The variable $v_i$ might occur in $\varphi$ in several places which
renders the straightforward substitution $\varphi \frac{t}{v_i}$ ``untidy''.
We ``tidy up'' $\varphi \frac{t}{v_i}$ by renaming variables. The assignment
of the variables of $t$ has to be modified accordingly.

So consider a formula $\varphi$, a variable $v_i $, a term $t$, and an
assignment $a$ with $\{ i | v_i \in \tmop{free} ( \varphi ) \} \cup \{ v_j \}
\subseteq \tmop{dom} ( a )$, where $v_j$ is a further variable thought to be a
bound for $t$ as in $\exists v_i \in v_j \varphi [ a ]$. Also assume that
$\varphi$ and $t$ have no common variable and do not contain $v_j $. Define
the {\tmstrong{tidy substitution}} $( \varphi \frac{t}{v_i} )^{\tmop{tidy}}$
of $v_i$ by $t$ into $\varphi$ as follows. If $v_i \nin \tmop{free} ( \varphi
)$ then let $( \varphi \frac{t}{v_i} )^{\tmop{tidy}} = \varphi$. Otherwise
ensure that $v_i$ is not a bound variable of $\varphi$ by possibly renaming
bound variables. Let $\varphi'$ be the renamed formula. Then rename all
occurances of $v_i$ in $\varphi'$ by {\tmem{pairwise different new}} variables
$w_0, \ldots, w_{k - 1} $, say. Obtain terms $t_0, \ldots, t_{k - 1}$ from the
given term $t$ by renaming all variables with {\tmem{new}} variables so that
for $i \neq j$ the terms $t_i$ and $t_j$ do not have common variables; call
$t_0, \ldots, t_{k - 1}$ {\tmstrong{copies}} of $t$. Now set
\[ \text{$( \varphi \frac{t}{v_i} )^{\tmop{tidy}} = \varphi' \frac{t_0 \ldots
   t_{k - 1}}{w_0 \ldots w_{k - 1}} $.} \]
The assignment $a$ has to be extended to an assignment $b \prec^{\ast} a$ in
line with the various renaming operations.

We define that the assignment $b$ {\tmstrong{adequately extends}} $a$
{\tmstrong{for the tidy substitution}} $( \varphi \frac{t}{v_i}
)^{\tmop{tidy}}$ {\tmstrong{bounded by}} $v_j$ if
\begin{enumeratealpha}
  \item $\forall l \in \tmop{dom} ( a ) \setminus \{ j \} b ( l ) = a ( l )$
  and $b ( j ) < a ( j )$;
  
  \item $\forall l \in \tmop{dom} ( b ) \setminus \tmop{dom} ( a ) b ( l ) < a
  ( j ) ;$
  
  \item if $v_l$ is a variable of $t$ and $v_{l'}$ and $v_{l''}$ are the
  renamings of $v_l$ in the copies $t_i$ and $t_j$ resp. then $b ( l' ) = b (
  l'' )$.
\end{enumeratealpha}
Note that $b$ satisfies $b \prec^{\ast} a$ by the substitution property Lemma
\ref{substition_property}.

With these preparations we can now carry out a recursive definition of the
bounded constructible truth function:
\begin{itemizeminus}
  \item $W ( a, \neg \varphi ) = 1$ iff $W ( a, \varphi ) = 0$;
  
  \item $W ( a, ( \varphi \wedge \psi ) ) = 1$ iff $W ( a, \varphi ) = 1$
  {\tmem{and}} $W ( a, \psi ) = 1$;
  
  \item $W ( a, \forall v_i \in v_j \varphi ) = 1$ iff {\tmem{for all}} terms
  $t$ and all assignments $b$ which are adequate for the tidy substitution $(
  \varphi \frac{t}{v_i} )^{\tmop{tidy}}$ bounded by $v_j$ \\
  holds $W ( b, ( \varphi \frac{t}{v_i} )^{\tmop{tidy}} ) = 1$;
  
  \item $W ( a, \exists v_i \in v_j \varphi ) = 1$ iff {\tmem{there is}} a
  term $t$ and an assignment $b$ which is adequate for the tidy substitution
  $( \varphi \frac{t}{v_i} )^{\tmop{tidy}}$ bounded by $v_j$ so that \\
  $W ( b, ( \varphi \frac{t}{v_i} )^{\tmop{tidy}} ) = 1$;
  
  \item $W ( a, v_i \in v_j ) = 1$ iff $a ( i ) < a ( j )$;
  
  \item $W ( a, v_i \in \{ v_j \in v_k | \varphi \} ) = 1$ iff $W ( a, \exists
  v_j \in v_k ( v_i = v_j \wedge \varphi ) ) = 1$;
  
  \item $W ( a, \{ v_i \in v_j | \varphi \} \in v_k ) = 1$ iff {\tmem{there
  is}} a term $t$ and an assignment $b$ which is adequate for the tidy
  substitution $( ( \{ v_i \in v_j | \varphi \} = v_l ) \frac{t}{v_l}
  )^{\tmop{tidy}}$ bounded by $v_k$ so that $W ( b, ( ( \{ v_i \in v_j |
  \varphi \} = v_l ) \frac{t}{v_l} )^{\tmop{tidy}} ) = 1$;
  
  \item $W ( a, \{ v_i \in v_j | \varphi \} \in \{ v_m \in v_n | \psi \} ) =
  1$ iff {\tmem{there is}} a term $t$ and an assignment $b$ which is adequate
  for the tidy substitution $( ( \{ v_i \in v_j | \varphi \} = v_m \wedge \psi
  ) \frac{t}{v_m} )^{\tmop{tidy}}$ bounded by $v_n$ so that $W ( b, ( ( \{ v_i
  \in v_j | \varphi \} = v_m \wedge \psi ) \frac{t}{v_m} )^{\tmop{tidy}} ) =
  1$;
  
  \item $W ( a, v_i = v_j ) = 1$ iff $a ( i ) = a ( j )$;
  
  \item $W ( a, v_i = \{ v_j \in v_k | \varphi \} ) = 1$ iff \\
  $W ( a, \forall v_l \in v_i \exists v_j \in v_k ( \varphi \wedge v_l = v_j
  ) ) = 1$ and \\
  $W ( a, \forall v_j \in v_k ( \varphi \rightarrow \exists v_l \in v_i v_l =
  v_j ) ) = 1$;
  
  \item $W ( a, \{ v_i \in v_j | \varphi \} = \{ v_m \in v_n | \psi \} ) = 1$
  iff $W ( a, ( \forall v_i \in v_j ( \varphi \rightarrow \exists v_m \in v_n
  ( v_i = v_m \wedge \psi ) ) ) = 1$ and $W ( a, \forall v_m \in v_n ( \psi
  \rightarrow \exists v_i \in v_j ( v_m = v_i \wedge \varphi ) ) ) ) = 1$.
\end{itemizeminus}
In all clauses the determination of $W ( a, \varphi )$ is reduced to values $W
( b, \psi )$. where the relevant arguments $( b, \psi )$ are of lesser
complexity than $( a, \varphi )$: either the assignments satisfy $b
\prec^{\ast} a$ or we have that $a = b$ and $\psi <_{\tmop{Form}} \varphi$.
Therefore $W$ has a recursive definition along the enumeration given by the
algorithm of the preceding section.

The programming techniques introduced above allow to incorporate the recursive
definition of $W$ into the enumeration algorithm for the class $A$ of
admissible pairs. This leads to our main results:

\begin{lemma}
  The bounded truth function $W$ for the constructible universe is ordinal
  computable.
\end{lemma}

\begin{theorem}
  \label{characterization}A set $x$ of ordinals is ordinal computable from a
  finite set of ordinal parameters if and only if it is an element of the
  constructible universe $L$. 
\end{theorem}

\begin{proof}
  Let $x \subseteq \tmop{Ord}$ be ordinal computable by the program $P$ from
  the finite set $\{ \alpha_0, \ldots, \alpha_{k - 1} \}$ of ordinal
  parameters: $P : \chi_{\text{$\{ \alpha_0, \ldots, \alpha_{k - 1} \}$}}
  \mapsto \chi_x $. By Lemma \ref{absolutness} c) the same computation can be
  carried out inside the inner model $L$:
  \[ ( L, \in ) \vDash P : \chi_{\text{$\{ \alpha_0, \ldots, \alpha_{k - 1}
     \}$}} \mapsto \chi_x . \]
  Hence $\chi_X \in L$ and $x \in L$.
  
  Conversely consider $x \in L$. Choose a tidy term $t$ and an assignment
  ($\alpha_0, \ldots, \alpha_{k - 1} ) \in \tmop{Ord}^{< \omega}$ such that $x
  = t [ ( \alpha_0, \ldots, \alpha_{k - 1} ) ]$. An ordinal $\beta$ can be
  represented as
  \begin{eqnarray*}
    \beta & = & \{ \alpha \in L_{\beta} | \alpha \text{ is an ordinal} \}\\
    & = & \{ v_{k + 1} \in v_k | v_{k + 1} \text{ is an ordinal} \} [ (
    \alpha_0, \ldots, \alpha_{k - 1}, \beta ) ] .
  \end{eqnarray*}
  Thus
  \begin{eqnarray*}
    \chi_x ( \beta ) = 1 & \leftrightarrow & \beta \in t [ ( \alpha_0,
    \ldots, \alpha_{k - 1} ) ]\\
    & \leftrightarrow & ( \{ v_{k + 1} \in v_k | v_{k + 1} \text{ is an
    ordinal} \} \in t ) [ ( \alpha_0, \ldots, \alpha_{k - 1}, \beta ) ]\\
    & \leftrightarrow & W ( ( \alpha_0, \ldots, \alpha_{k - 1}, \beta ),
    \{ v_{k + 1} \in v_k | v_{k + 1} \text{ is an ordinal} \} \in t ) = 1.
  \end{eqnarray*}
  Using the enumeration algorithm for the truth function $W$ one can turn this
  equivalence into an ordinal algorithm which sends $\chi_{\{ \alpha_0,
  \ldots, \alpha_{k - 1} \}}$ to $\chi_x$. Hence $x$ is ordinal computable
  from the parameters $\alpha_0, \ldots, \alpha_{k - 1} \in \tmop{Ord}$.
\end{proof}

\section{The Generalized Continuum Hypothesis in $L$}

Ordinal computability allows to reprove some basic facts about the
constructible universe $L$. The analogue of the axiom of constructibility, $V
= L$, is the statement that {\tmem{every}} set of of ordinals is ordinal
computable from a finite set of ordinals.

\begin{theorem}
  The constructible model $( L, \in )$ satisfies that every set of ordinals is
  ordinal computable from a finite set of ordinals.
\end{theorem}

\begin{proof}
  Let $x \in L$, $x \subseteq \tmop{Ord}$. By Theorem \ref{characterization},
  take a program $P$ and a finite set $\{ \alpha_0, \ldots, \alpha_{k - 1} \}$
  of ordinal parameters such that $P : \chi_{\text{$\{ \alpha_0, \ldots,
  \alpha_{k - 1} \}$}} \mapsto \chi_x $. By Lemma \ref{absolutness} c) the
  same computation can be carried out inside the inner model $L$:
  \[ ( L, \in ) \vDash P : \chi_{\text{$\{ \alpha_0, \ldots, \alpha_{k - 1}
     \}$}} \mapsto \chi_x . \]
  So in $L$, $x$ is ordinal computable from the set $\{ \alpha_0, \ldots,
  \alpha_{k - 1} \}$.
\end{proof}

The following therem is proved by a condensation argument for ordinal
computations which is a simple analogue of the usual condensation arguments
for the constructible hierarchy.

\begin{theorem}
  Assume that every set of ordinals is ordinal computable from a finite set of
  ordinals. Then:
  \begin{enumeratealpha}
    \item Let $\kappa \geqslant \omega$ be an infinite ordinal and $x
    \subseteq \kappa$. Then there are ordinals $\alpha_0, \ldots, \alpha_{k -
    1} < \kappa^+$ such that $x$ is ordinal computable from the set\\ $\{
    \alpha_0, \ldots, \alpha_{k - 1} \}$.
    
    \item Let $\kappa \geqslant \omega$ be infinite. Then $\tmop{card}
    (\mathcal{P}( \kappa ) ) = \kappa^+$.
    
    \item The generalized continuum hypothesis $\tmop{GCH}$ holds.
  \end{enumeratealpha}
\end{theorem}

\begin{proof}
  a) Take a program $P$ and a finite set $\{ \alpha'_0, \ldots, \alpha'_{k -
  1} \}$ of ordinal parameters such that $P : \chi_{\text{$\{ \alpha'_0,
  \ldots, \alpha'_{k - 1} \}$}} \mapsto \chi_x $; let $\theta$ be the length
  of this ordinal computation. Take a transitive $\tmop{ZF}^-$-model $( M, \in
  )$ such that $\alpha'_0, \ldots, \alpha'_{k - 1}, \theta, \kappa, x \in M$.
  By Lemma \ref{absolutness} a), $( M, \in )$ also satisfies that $P :
  \chi_{\text{$\{ \alpha'_0, \ldots, \alpha'_{k - 1} \}$}} \mapsto \chi_x $.
  The downward {\tmname{Löwenheim-Skolem}} theorem and the
  {\tmname{Mostowski}} isomorphism theorem yield an elementary embedding
  \[ \pi : ( \bar{M}, \in ) \rightarrow ( M, \in ) \]
  such that $\bar{M}$ is transitive, $\tmop{card} ( \bar{M} ) = \kappa$ and
  $\{ \alpha'_0, \ldots, \alpha'_{k - 1}, \theta, \kappa, x \} \cup \kappa
  \subseteq \pi'' \bar{M}$. Let $\pi ( \alpha_0 ) = \alpha_0', \ldots, \pi (
  \alpha_{k - 1}' ) = \alpha_{k - 1} $. Then $\alpha_0, \ldots, \alpha_{k - 1}
  < \kappa^+$ since $\tmop{card} ( \bar{M} ) < \kappa^+$. Observe that $\pi (
  x ) = x$. Since $\pi$ is elementary $( \bar{M}, \in )$ satisfies that $P :
  \chi_{\text{$\{ \alpha_0, \ldots, \alpha_{k - 1}' \}$}} \mapsto \chi_x $. By
  Lemma \ref{absolutness} d), $P : \chi_{\text{$\{ \alpha_0, \ldots, \alpha_{k
  - 1}' \}$}} \mapsto \chi_x$ in $V$. Thus $x$ is ordinal computable from the
  set $\{ \alpha_0, \ldots, \alpha_{k - 1} \}$ as required.\\
  b) follows from a) since there are a countable many programs and $\kappa^+$
  many finite sets of ordinals $< \kappa^+$.\\
  c) is immediate from b)
\end{proof}

These two theorems immediately imply {\tmname{Gödel}}'s result:

\begin{theorem}
  $( L, \in ) \vDash \tmop{GCH}$.
\end{theorem}

Other condensation arguments like the proof of the combinatorial principle
$\diamondsuit$ in $L$ can also be translated into the setting of ordinal
computability in a straightforward way. It remains to be seen whether
arguments involving {\tmname{Jensen}}'s fine structure theory of the
constructible hierarchy {\cite{Jensen}} can be carried out with ordinal
computability. One would hope that the simple concept of ordinal computation
allows clear proofs of principles like $\Box$ and morasses without
definability complications.

\end{document}